\documentclass[11pt]{article}
\usepackage{amsmath}
\usepackage{amssymb}
\usepackage{amsthm}
\usepackage{amscd}
\usepackage{epsfig}

\newcommand{\remove}[1]{}
\newtheorem {Theorem}    {Theorem}[section]
\newtheorem {Definition} {Definition} [section]

\newtheorem {Lemma}      [Theorem]    {Lemma}

\newtheorem {Proposition}[Theorem]    {Proposition}
\newtheorem {Claim}      [Theorem]    {Claim}

\newcommand {\Heads}[1]   {\smallskip\pagebreak[1]\noindent{\bf
#1{\hskip 0.2cm} }}
\newcommand {\Head}[1]    {\Heads{#1:}}

\newcommand {\proofs}     {\proof}
\newcommand {\prooft}[1]  {\Head{Proof #1}\nopagebreak[2]}

\newcommand  {\QED}    {\def\qedsymbol{$\blacksquare$}\qed}

\newcommand{\ignore}[1]{}

 \def\eps{\epsilon} 
\def\P{{\bf{P}}}
\def\V{{\bf{V}}}

\def\R{\hbox{I\kern-.2em\hbox{R}}}
\def\A{{\cal{A}}}

\def\F{{\hat{F}}}

\def\B{{\cal{B}}}
\def\|{\, | \, }
\def\v0{{\bf 0}}

\def\0{\hat{0}}
\def\1{\hat{1}}

\def\path{{\tt path}}

\def\phi{\varphi}

\def\be{\begin{equation}}
\def\ee{\end{equation}}

\def\L{{\cal L}}
\def\d{{\hat{d}}}
\def\B{{\hat{B}}}
\def\net{{\cal N}}
\def\part{{\cal P}}
\def\vert{{\cal V}}
\def\edges{{\cal E}}

\title{Distorted metrics on trees and phylogenetic forests}

\author{{\sc Elchanan Mossel} \thanks{{\tt mossel@stat.berkeley.edu} 
Miller fellow, CS and Statistics, U.C. Berkeley}}
\date{}
\begin{document}
\maketitle

\begin{abstract}
We study distorted metrics on binary trees in the context 
of phylogenetic reconstruction. 
Given a binary tree $T$ on $n$ leaves with a path metric $d$, 
consider the pairwise distances $\{d(u,v)\}$ between
leaves. It is well known that these determine the tree and the $d$
length of all edges. Here we consider distortions $\d$ of $d$ such that 
for all leaves $u$ and $v$ it holds that $|d(u,v) - \d(u,v)| < f/2$ if
either $d(u,v) < M$ or $\d(u,v) < M$, where $d$ satisfies $f \leq d(e)
\leq g$ for all edges $e$. 
Given such distortions we show how to reconstruct in polynomial time a forest
$T_1,\ldots,T_{\alpha}$ such that the true tree $T$ may be obtained from
that forest by adding $\alpha-1$ edges and 
$\alpha-1 \leq 2^{-\Omega(M/g)} n$.

Metric distortions arise naturally in phylogeny, 
where $d(u,v)$ is defined by the log-det 
of a covariance matrix associated with $u$ and $v$. 
When $u$ and $v$ are ``far'', the entries of the covariance matrix are
small and therefore $\d(u,v)$, which is defined by log-det of an
associated empirical-correlation matrix may 
be a bad estimate of $d(u,v)$ even if the correlation matrix is ``close'' to 
the covariance matrix. 

Our metric results are used in order to show how to reconstruct
phylogenetic forests with small number of trees from sequences of length 
logarithmic in the size of the tree. Our method also yields 
an independent proof that phylogenetic
trees can be reconstructed in polynomial time from sequences of
polynomial length under the standard assumptions in phylogeny. Both the metric
result and its applications to phylogeny are almost tight. 
\end{abstract}

%%%%%%%%%%%%%%%%%%%%%%%%%%%%%%%%%%%%%%%%%%%%%%%%%%%%%%%%%%%%%%%%%%%%%
%%%%%%%%%%%%%%%%%%%%%%%%%%%%%%%%%%%%%%%%%%%%%%%%%%%%%%%%%%%%%%%%%%%%%
\section{Introduction}
%%%%%%%%%%%%%%%%%%%%%%%%%%%%%%%%%%%%%%%%%%%%%%%%%%%%%%%%%%%%%%%%%%%%%
%%%%%%%%%%%%%%%%%%%%%%%%%%%%%%%%%%%%%%%%%%%%%%%%%%%%%%%%%%%%%%%%%%%%%

Reconstructing 
phylogenies have been a scientific challenge for the last 50 years. 
We refer the reader to \cite{Felsenstein:04} or \cite{SempleSteel:03} 
for general and
mathematical background. 
The standard setting in phylogeny is of trees where the leaves are
labeled by {\em taxa} or {\em species}. 
Given aligned sequences at the leaves, we define character $i$ to be 
the collection of letters at position $i$ for all the species.
Under the i.i.d. assumption,  
the characters are independent samples from the
evolutionary process on the tree.

The theoretical foundations of most methods  
used in phylogeny are unsatisfactory. 
Under the standard i.i.d. model, {\em Parsimony} is
not consistent \cite{Cavender:78,Felsenstein:78} 
and is NP hard to compute \cite{DaJoSa:86,FouldsGraham:82,GrahamFoulds:82}. 
The computational complexity of finding the {\em Maximum likelihood} 
tree is not known and the best bounds on the the amount of data 
needed are {\em exponential} in the number of taxa \cite{SteelSzekely:03}.

Computational complexity and information theory considerations have not played
an important role in phylogeny in the past as biologists were
mostly interested in reconstructing trees on a small number 
(typically at most a few hundred) species.
However, one of the major goals of systematic biology in the coming decade is 
to reconstruct phylogenies on millions of species. It is clear
that for such numbers, it is crucial to apply algorithms with low
computational complexity. Similarly, algorithms should use information
efficiently.

In \cite{ErStSzWa:99a} the authors developed the first
reconstruction algorithm satisfying two important properties:
\begin{itemize}
\item
Given number of characters that is polynomial in the number of taxa,
the algorithm finds the true tree with high probability.
\item
The running time of the algorithm is polynomial in the number of taxa.
\end{itemize}
Variants of these method and generalizations from two states models to
general models appeared in \cite{ErStSzWa:99b}. In \cite{CrGoGo:03} the authors
discuss a closely related problem of learning a
phylogenetic tree (in the PAC setting). They developed a PAC learning
algorithm for the two state model. The problem of PAC-learning general-state
model in polynomial time is still open. See also \cite{FarachKannan:99} for an
earlier result on learning phylogenies.

The method developed in \cite{ErStSzWa:99a} is a {\em distance}
method. Such methods were commonly used in phylogeny before, but
\cite{ErStSzWa:99a} is the result where a distance method yields a provably 
good performance. Distance methods are based on defining a path metric
on the tree based on the evolution model. Then the distance between
leaves of the tree is approximated by some distance between the
corresponding sequences at the leaves. In this sense all distance
methods in phylogeny 
may be view as reconstruction methods from distorted metrics.
%While it is clear that any reconstruction method will run in time
%at least polynomial in the number of taxa (since reading all the taxa takes a
%polynomial time), it is far from clear that the number of characters
%needed is indeed polynomial in the number of taxa. 

Given the existence of a polynomial time reconstruction algorithms, 
the next problem is optimizing
the sampling complexity.
The number of characters needed ($=$ the length of sequences) 
is of great practical importance, as
this number is bounded by the underlying biology. It is therefore
desirable to minimize this number. 

Since there are $2^{\Theta(n \log n)}$ binary trees on $n$ leaves,
an easy counting argument yields that the number of characters needed
is at least logarithmic in the number of taxa. Thus we are led to the
following natural problem:
is the length of the sequences needed logarithmic in $n$ as
or is it polynomial?

%As an example note that
%the value of the polynomial $p(n) = n^2$ for $n=10^6$ is $10^{12}$,
%while $\log( 10^6 ) \leq 20$. It is highly unlikely that we could find
%aligned sequences of length $10^{12}$ for $10^6$ species. Sequences of
%length $20$ are much more likely to be found. 

In \cite{Mossel:03a} we showed that for a restricted family of models, it
is possible to reconstruct phylogenies from a logarithmic number of
characters, if the mutation rates are low (bounded above by some
constant). We also showed that a 
polynomial number of characters is needed if the
mutation rates are high (bounded below by some constant). 

We later \cite{Mossel:03b} (see also \cite{SoberSteel:02}) 
generalized the polynomial lower bound for high mutation
rates to a large family of models. In \cite{MosselSteel:03} we analyze 
another model where logarithmic reconstruction is achievable for low
mutation rates. 

The phase transition discussed above is of crucial interest if we wish
to reconstruct all the tree. However, in some cases, a more modest
objective is posed: reconstruct a ``large portion'' of the 
tree. Practitioners (this was kindly noted to me by J. Felsenstein 
(2001, private communication) and J. Kim (2003, private communication)) 
have noticed that this problem seems to be much
easier than the problem of reconstructing the complete tree. 

\begin{figure}[h] \begin{center} \label{add_edges}
\input{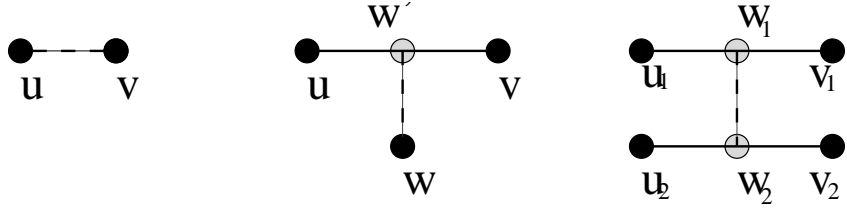}
\caption{The ``edge-adding'' operation}
\end{center}
 \end{figure}

In this paper we prove that this is indeed the case. 
\begin{Definition} \label{def:edge_add}
We define the
operation of {\em edge adding} to a forest as one of the following
\begin{itemize}
\item
Add an edge $(u,v)$ connecting two isolated leaves $u$ and $v$.
\item
Given an edge $(u,v)$ of the forest and an isolated leaf $w$,
replace the edge $(u,v)$ by the edges $(u,w'),(w',v)$ and $(w',w)$
where $w'$ is a new vertex.
\item
Replace the two edges $(u_1,v_1)$ and $(u_2,v_2)$ of the forest by  
$(u_1,w_1),(w_1,v_1),(u_2,w_2),(w_2,v_2)$ and $(w_1,w_2)$ where
$w_1,w_2$ are new vertices.
\end{itemize}
See figure \ref{add_edges}.
\end{Definition}

We show that under the standard assumptions in phylogeny, 
for a tree $T$ on $n$ leaves and for all $\delta > 0$ and $\gamma > 0$
we can reconstruct from $\gamma^{-O_{\delta}(1)} \log n$ characters 
a forest $T_1,\ldots,T_{\alpha}$ such that 
${\alpha} \leq 1 + \gamma n$ 
and that $T$ may be obtained from the forest $T_1,\ldots,T_{\alpha}$
by adding at most $\alpha-1$ edges. 
The reconstruction is performed in a polynomial time and with error
bounded by $\delta$.

Note that taking $\gamma$ to be a small constant we obtain that ``most''
edges of the tree can be reconstructed from $O(\log n)$ characters. 
Taking $\gamma = 1/n$, we obtain that the full tree may be
reconstructed from a polynomial number of characters. Thus obtaining an
independent proof of the results of \cite{ErStSzWa:99a,ErStSzWa:99b,CrGoGo:03}. 

Our results indicate what level of refinement is achievable in
reconstructing a phylogenetic tree given a certain amount of data. 
%Indeed, below we will show that for general mutation models there
%results are tight up to the constant appearing in the $O(1)$ above.
We believe that our techniques may also play an important role 
in reconstructing the complete phylogenetic tree in low mutation rates
from logarithmic number of characters 
as it allows a very clean divide and conquer approach (see \cite{Mossel:03a}). 

The main ingredient of the proof is a metric theorem that
can be stated roughly as follows. Let $T$ be a binary tree and $d$ a path
metric on $T$ such that $f \leq d(e) \leq g$ for all edges $e$. Let
$\L(T)$ be the set of leaves of $T$ and $\d : \L(T) \times \L(T) \to
\R$ a distortion of $d$ that satisfies $|d(u,v) - \d(u,v)| < f/2$ if 
$\d(u,v) < M$ or $d(u,v) < M$ ($\d$ typically does not correspond to a 
path metric on $T$).
We show that we can partition $\L(T)$ into
$\alpha$ sets $L_1,\ldots,L_{\alpha}$, where $\alpha = n
2^{-\Omega(M/g)}$. For $1 \leq \beta \leq \alpha$ we reconstruct a
tree such that the leaves of $T_{\beta}$ are 
$L_{\beta}$ and the tree $T$ may be obtained from the forest $(T_{\beta})_{\beta \leq \alpha}$ 
by adding at most $\alpha-1$ edges. 

It is easy to see that the metric theorem is tight up to the 
constant in the $\Omega$ by 
considering an $r$-level $3$-regular tree, where the length of all
edges is exactly $g$ and $\d(u,v) = \infty$ if $d(u,v) > M$.
Similar tightness results hold for the number of trees in the forest in the phylogenetic reconstruction. The proof is more complicated and follows ideas from 
\cite{Mossel:03b} on lower bounds on 
the  sampling complexity of phylogenies. The proof is omitted in this extended abstract.
Our metric result may be of independent interest to other problems
where path metrics on trees are considered.

We now give a high level sketch of the different sections of the paper. 
\begin{itemize}
\item
The formal definition of the model and the statement of the main results 
are given in Section \ref{sec:def}, where we also discuss how the metric
result implies the result in phylogeny. 

The reduction from the phylogenetic model to distorted metrics is
given in Proposition \ref{prop:large_deviations} which is an easy 
reformulation of a large deviation result from \cite{ErStSzWa:99b}. 
This proposition gives $M$ as a function of the number of samples. 
The distortion error is assumed to be at most $\eps=f/2$. 
It is known that Proposition \ref{prop:large_deviations} is
essentially tight. In other words, distances larger than $M$ are
likely to be computed with error larger than $f/2$. 
\item
Given the distortion, it is easy to construct for  
each leaf $v$, a tree $T_v$ on the set of leaves in the $(M-7f)/6$ (metric)-neighborhood of $v$.
This is done using standard techniques in Phylogeny.
\item
The collection of trees $\{T_v\}$ is {\em not} the forest we are
looking for. First, the trees in this forest are {\em not} disjoint.
Second, there may be many trees in this collection (in fact as many as 
$n$ different trees). 
The main task of the paper is to ``glue'' these trees to form {\em
  edge-disjoint} forest. Then we can bound the number of trees in the forest.
\item
The notion of edge-disjoint trees is studied in Section
\ref{sec:edge_disjoint}. The results of this section imply that if a forest
of edge-disjoint trees is a 
refinement of the collection $\{T_v\}$ above then the
size of the forest is $1 + n2^{-\Omega(M/g)}$.
\item
The ``glueing'' algorithm is given in Section \ref{sec:super_tree}. 
\item
The final metric result is stated in Section \ref{sec:dist}.
\end{itemize}

{\bf Acknowledgments:}
The idea that reconstructing a forest should be ``easy'' was conceived during
a talk by Junhyong Kim at the kickoff meeting of the Cipres project. J. Kim 
said that it seems like most edges of phylogenies are 
easy to reconstruct. I thank him and Tandy Warnow for encouragement 
to work on this problem.

%%%%%%%%%%%%%%%%%%%%%%%%%%%%%%%%%%%%%%%%%%%%%%%%%%%%%%%%%%%%%%%%%%%%%
%%%%%%%%%%%%%%%%%%%%%%%%%%%%%%%%%%%%%%%%%%%%%%%%%%%%%%%%%%%%%%%%%%%%%
\section{Definitions and main results} \label{sec:def}
%%%%%%%%%%%%%%%%%%%%%%%%%%%%%%%%%%%%%%%%%%%%%%%%%%%%%%%%%%%%%%%%%%%%%
%%%%%%%%%%%%%%%%%%%%%%%%%%%%%%%%%%%%%%%%%%%%%%%%%%%%%%%%%%%%%%%%%%%%%
Let $T$ be a tree. Write $\vert(T)$ for the nodes of $T$,
$\edges(T)$ for the edges of $T$ and $\L(T)$ for the leaves of $T$.
If the tree is rooted, then we denote by $\rho(T)$ the root of $T$.  
Unless stated otherwise, all trees are assumed to be {\em binary} 
(all internal degrees are $3$) 
and it is further assumed that $\L(T)$ is labeled. 

Let $T$ be a tree equipped with a path metric 
$d : \edges(T) \to \R_{+}$. $d$ will also denote the induced metric 
on $\vert(T)$:
\begin{equation} \label{eq:path_metric}
d(v,w) = \sum \{d(e) : e \in \path(v,w)\},
\end{equation}
for all $v,w \in \vert(T)$.

We will further assume below that the length of all edges 
is bounded between $f$ and $g$ for all $e \in E$. In other words, for
all $e \in \edges(T)$,
\begin{equation} \label{eq:fg}
f \leq d(e) \leq g.
\end{equation}

In applications to phylogeny we are typically given a distortion
$\d : \L(T) \times \L(T) \to \R_{+}$ of $d$. We define an $(\eps,M)$
distortion as follows. 
%(typically $\d$ does not correspond to any path
%metric on the tree $T$).

\begin{Definition} \label{def:distortion2}
Given a tree $T$ equipped with a metric $d$, and two positive
numbers $0 < \eps < M$, we say that 
$\d : \L(T) \times \L(T) \to \R_{+} \cup \{\infty\}$ is an 
$(\eps,M)$ distortion of $d$ if 
\begin{itemize}
\item
$\d(u,v) = \d(v,u)$ for all $u$ and $v$ in $\L(T)$; i.e., $\d$ is symmetric.
\item
If $\d(u,v) = \infty$, then $d(u,v) > M$. 
\item
If $\d(u,v) < \infty$, then  
$|\d(u,v) - d(u,v)| < \eps$. 
\end{itemize}
\end{Definition}

It is well known that $d : \L(T) \times \L(T) \to \R_{+}$ determines
the underlying tree $T$ and the metric on the
edges $d : \edges(T) \to \R_{+}$. Moreover, there exists a polynomial 
time algorithm to reconstruct $T$. 
Similarly, we may recover $T$ 
from any $(\eps,\infty)$ distortion of $d$ if $\eps < f/2$. 
Moreover, in this case, we may also recover a function $\d : \edges(T) \to \R$
satisfying $|\d(e) - d(e)| < 2\eps$ 
(\cite{Buneman:71,Gusfield:91}, see e.g. \cite[Chapter 7]{SempleSteel:03}). 

In our main result we show that given an $(\eps,M)$ distortion of $d$, 
we may recover many of the edges of $T$ and a 
good approximation of the $d$ length of those edges. 

\begin{Theorem} \label{thm:part}
Let $T=(V,E)$ be a binary tree equipped with a metric $d$ satisfying 
(\ref{eq:fg}). Let $n = |\L(T)|$, so that $|\vert(T)| = 2 n - 2$. 
Let $\d$ be an $(\eps,M)$ distortion of $d$ and suppose
that $\eps < f/2$ and that $M > 7 \eps$. 
Then $\d$ determines a partition $\part$ of $\L(T)$ into sets 
$L_1,\ldots,L_{\alpha}$ and a forest $T_1,\ldots,T_{\alpha}$ such that 
$L_{\beta} = \vert(T_{\beta}) \cap \L(T)$ for all $\beta$ and  
\begin{itemize}
\item
The tree $T$ may be obtained from the forest 
$T_1,\ldots,T_{\alpha}$ by adding at most  
$\alpha - 1$ edges. 
\item
The number of trees $\alpha$ in the forest is 
at most $\lfloor 1 + \frac{60 n}{\sqrt{2}} 2^{-\frac{M-\eps}{2g}} \rfloor$. 
\end{itemize}

Moreover, 
\begin{itemize}
\item
the partition $(L_{\beta})_{\beta \leq \alpha}$,
\item 
the trees $(T_{\beta})_{\beta \leq \alpha}$ and 
\item
a function 
$\d : \cup_{\beta \leq \alpha} \edges(T_{\beta}) : \to \R_{+}$
satisfying $|\d(e)-d(e)| < 2 \eps$,
\end{itemize}
can be all computed from $\d$ in time polynomial on $n$.
\end{Theorem}

{\bf Mutation models and distances.} 
When reconstructing phylogenies, the data is given as
sequences at the labeled leaves $\L(T)$ and the tree $T$ is unknown.
Usually the mutation model is defined on a rooted tree while 
the goal is to reconstruct un-rooted trees (in many
models there is no way to distinguish a root). 

We let $\A$ denote the alphabet in which information is encoded. For example, 
$\A = \{A,C,G,T\}$ for DNA sequences, $\A = \{20 \mbox{ amino acides}
\}$ for proteins
and $\A = \{0,1\}$ for purine-pyrmidine sequences. To define the
mutation model we assume that all the edges are directed away from the
root and for each edge $e \in \edges(T)$ let $M(e)$ be the mutation
matrix corresponding to the edge $e$. $M(e)$ is an $|\A| \times |\A|$ 
stochastic matrix. The $(i,j)$'th entry of $M(e)$ is the probability that state
$i$ will mutate to state $j$ along edge $e$. It is assumed that each
character evolves down the tree as a Markov-chain on the tree, where
$M(e)$ is the transition matrix for edge $e$. The root letter is
chosen from some fixed distribution $\pi$. It is assumed that the
characters evolve in an i.i.d. manner - they all come 
from the same distribution and each one is independent from all the others.

Two popular examples are the CFN model where 
$M(e) = \left( \begin{matrix} 1 - \theta(e) & \theta(e) \\ 
                                \theta(e) & 1 - \theta(e) \end{matrix}
                              \right)$
and the Jukes-Cantor models where 
\[
M(e) = \left( \begin{matrix} 
1 - 3\theta(e) & \theta(e) & \theta(e) & \theta(e) \\
\theta(e) & 1 - 3\theta(e) & \theta(e) & \theta(e) \\
\theta(e) & \theta(e) & 1 - 3\theta(e) & \theta(e) \\
\theta(e) & \theta(e) & \theta(e) & 1 - 3\theta(e) \\
\end{matrix} \right).
\]

It turns out that under mild assumptions on the matrices $M(e)$ and
the evolution model - the {\em log-det} distance defines a path metric
on the tree \cite{Steel:94}. We summarize the basic properties of this distance if the
following proposition. The proof for the CFN and the Jukes-Cantor
model have appeared independently several times. The general case 
follows from a large deviation estimate in \cite{ErStSzWa:99b} and is proven
in the appendix.
\begin{Proposition} \label{prop:large_deviations}
Assume that the matrices $M(e)$ satisfy that 
$e^{-2g'} < \det(M(e)) < e^{-2f}$ for all 
$e \in \edges(T)$ and that for all nodes
$v \in \vert(T)$ and all letters $a \in \A$, the probability that the
letter at $v$ is $a$ is at least $\pi_{\min} > 0$. Let $\eps > 0$.

For every two vertices $u,v \in \L(T)$ and $a,b \in \A$, 
let $F_{a,b}$ be the probability that node $u$ has letter $a$ and node
$v$ has letter $b$. Let $\F_{a,b}$ be the empirical distribution that
node $u$ has letter $a$ and node $v$ has letter $b$.  
Let $d(u,v) = -\log \det (F_{i,j})$ and 
$\d(u,v) = -\log \det (\F_{i,j})$ if $-\log \det (\F_{i,j}) \leq M+\eps$,
and 
$\d(u,v) = \infty$ if $-\log \det (\F_{i,j}) > M+\eps$ or 
$\det (\F_{i,j}) \leq 0$.
Then 
\begin{itemize}
\item
$d(u,v)$ is a path metric on the tree satisfying  
$g \geq d(e) \geq f$ for all edges $e$ of the tree (where $g$ depends
on $g'$ and $\pi_{\min}$). 
\item
There exists a constant $c$ such that for all $r > 2$ if the number of sample 
satisfies  
\begin{equation} \label{eq:2nd_dev}
k \geq \frac{c r}{(1 - e^{-2\eps})^2} e^{2M+2\eps} \log n, 
\end{equation}
then with probability at least 
$1 - n^{2 - r}$ it holds that 
$\d$ is an $(\eps,M)$ distortion of $d$. 
\end{itemize}
\end{Proposition}

The proposition may be used with different values of the parameters 
$k$ and $M$. Fix $\eps < f/2$. Taking $M$ to be the diameter of the tree, 
it gives that all empirical distances are within $\eps$ 
of the true distances once $k$ is exponential in $M$. Since $M$ may be as 
large as $\Omega(g n)$, this gives sampling complexity 
$k$ which is exponential in $n$. 
Taking $M = 100g$, say, would give an $(\eps,100g)$ distortion from 
$k = O(e^{100 g} \log n) = O(\log n)$ samples. In sequel we will use  
$M$ ranging between $M = O(g)$ for a constant $c$ to $k = O(g \log n)$ 
(in particular, typically we will only have 
a fraction of the distances within $\eps$ of their true value). 

Combining Proposition \ref{prop:large_deviations} and Theorem
\ref{thm:part} we obtain
\begin{Theorem} \label{thm:phyl}
Consider a binary phylogenetic tree $T$,
where the log-det distance associated with the mutation matrices 
$M(e)$ satisfy the conditions of Proposition
\ref{prop:large_deviations}. Then given $k$ satisfying
(\ref{eq:2nd_dev}), with $M > 7 \eps$, and $\eps < f/2$ 
we can with probability at least $1 - n^{2-r}$ recover a partition 
$\part$ of $\L(T)$ into sets 
$L_1,\ldots,L_{\alpha}$ and a forest $T_1,\ldots,T_{\alpha}$ such that 
$L_{\beta} = \vert(T_{\beta}) \cap \L(T)$ for all $\beta$ and 
\begin{itemize}
\item
$T_1,\ldots,T_{\alpha}$ is a forest that may be 
obtained from $T$ by removing $\alpha - 1$ edges. 
\item
The number of trees $\alpha$ in the forest is 
at most $\lfloor 1 + \frac{60 n}{\sqrt{2}} 2^{-\frac{M-\eps}{2g}} \rfloor$. 
\end{itemize} 
Moreover, we can recover a function 
$\d : \cup_{\beta \leq \alpha} \edges(T_{\beta}) : \to \R_{+}$
satisfying $|\d(e)-d(e)| < 2 \eps$ for all $e$.
\end{Theorem}

Note that taking $M = O(g)$ to be a large constant and $\eps = f/2$ 
proves that most edges of the tree can be recovered from $O(\log n)$
characters. Similarly, taking $M = O(g \log n)$ and $\eps = f/2$, we
see that we can recover the underlying tree from $k = n^{O(g)}$ characters, 
thus obtaining an independent proof of the results of \cite{ErStSzWa:99a,ErStSzWa:99b,CrGoGo:03}.

%%%%%%%%%%%%%%%%%%%%%%%%%%%%%%%%%%%%%%%%%%%%%%%%%%%%%%%%%%%%%%%%%%%%%
%%%%%%%%%%%%%%%%%%%%%%%%%%%%%%%%%%%%%%%%%%%%%%%%%%%%%%%%%%%%%%%%%%%%%
\section{Edge disjoint trees} \label{sec:edge_disjoint}
%%%%%%%%%%%%%%%%%%%%%%%%%%%%%%%%%%%%%%%%%%%%%%%%%%%%%%%%%%%%%%%%%%%%%
%%%%%%%%%%%%%%%%%%%%%%%%%%%%%%%%%%%%%%%%%%%%%%%%%%%%%%%%%%%%%%%%%%%%%
{\em Edge disjoint trees} and {\em edge sharing trees} 
will play a crucial role below. In this section we define these notions and
discuss some of their basic properties. 
% VERSIONS %
%Most of the proofs in this section are omitted and can be found in the appendix. 
% END VERSIONS %

We let $T$ be a binary tree with vertices $\vert(T)$
and edges $\edges(T)$. 
We let $\L(T) \subset \vert(T)$ be the set of leaves of $T$ and 
$n = |\L(T)|$ the size of this set.
We write $\path_T(x,y)$ for the path (sequence of edges) connecting $x$
to $y$ in $T$. We will sometime omit the subscript $T$ 
and write $\path(x,y)$. We write $\ell(x,y)$ or $\ell_T(x,y)$ 
for the number of edges in the path connecting $x$ and $y$. For two 
sets $A,B \subset \V(T)$ we write 
$\ell_T(A,B) = \min\{ \ell(x,y) : x \in A, y \in B\}$.

Removing an edge $e$ from a tree $T$ results in obtaining two trees 
$T_1$ and $T_2$. The {\em split defined by $e$} 
is the partition $\{\L(T) \cap \vert(T_1), \L(T) \cap \vert(T_2)\}$ 
of $\L(T)$. We denote by $\Sigma(T)$ the collection of $\L(T)$ splits
defined by all edges of $T$. It is well know that $\Sigma(T)$
determines $T$ (see e.g. \cite[Chapter 3]{SempleSteel:03}). We denote the split
$\{A,B\}$ by $A|B$.

\begin{Definition} \label{def:rest}
Let $T$ be a binary tree and $L \subset \L(T)$ a set of leaves. 
The {\em restriction} of $T$ to $L$ is defined as follows. This is the tree
whose leave set is $L$ and whose splits are defined by
\[
\Sigma(T|L) := \{A|A' : A = B \cap L, A' = B' \cap L \mbox{ and } B | B' \in
\Sigma(T)\}.
\]
The restriction of $T$ to $L$ is denoted by $T|L$.

Given two sets $L_1,L_2 \subset \L(T)$, we say that the trees 
$T|L_1, T|L_2$ are {\em edge disjoint} if 
\[
\path_T(u_1,v_1) \cap \path_T(u_2,v_2) = \emptyset,
\]
for all $u_1,v_1 \in L_1$ and $u_2,v_2 \in L_2$.
We say that $T|L_1,T|L_2$ are edge-sharing if they are not edge disjoint.
\end{Definition}

% VERSIONS - SHORT
%The following easy lemma whose proof is omitted 
%is useful as it shows that edge disjointness doesn't depend on the ambient tree.
%ENDSHROT

%VERSIONS - FULL
The following easy lemma is useful as it shows that edge disjointness
does not depend on the ambient tree.
%ENDFULL
\begin{Lemma} \label{lem:prop_disj}
Let $L_1 \cup L_2 \subset L' \subset \L(T)$. Then $T|L_1$ and $T|L_2$ 
are edge disjoint if and only if 
$\path_{T | L'}(u_1,v_1) \cap \path_{T | L'}(u_2,v_2) = \emptyset$, for all 
$u_1, v_1 \in L_1$ and $u_2, v_2 \in L_2$.

In particular, $L_1$ and $L_2$ are edge disjoint if and only if 
$\path_{T | L_1 \cup L_2}(u_1,v_1) \cap 
 \path_{T | L_1 \cup L_2}(u_2,v_2) = \emptyset$, 
for all $u_1, v_1 \in L_1$ and $u_2, v_2 \in L_2$.
\end{Lemma}

%%%%%%%%%%%%% VERSIONS SHORT
%\ignore{
%%%%%%%%%%%%% END SHORT
\proofs
The second statement follows immediately from the first one by letting 
$L' = L_1 \cup L_2$. 

For the first statement, note that $\path_T(u_1,v_1)$ is obtained from 
$\path_{T | L'}(u_1,v_1)$ by replacing each edge $(x,y) \in
\path_{T | L'} (u_1,v_1)$ by a sequence of edges 
$(x = x_1,x_2),\ldots,(x_j,x_{j+1} = y)$. Moreover, the sequence 
$(x = x_1,x_2),\ldots,(x_j,x_{j+1} = y)$ depends on the edge $(x,y)$
only and each edge $(x_i,x_{i+1})$ of $T$ appears in the sequence of at most  
one edge $(x,y)$ of $T | L'$.

It now follows that 
$\path_{T | L'}(u_1,v_1) \cap 
 \path_{T | L'}(u_2,v_2) = \emptyset$, 
for all $u_1, v_1 \in L_1$ and $u_2, v_2 \in L_2$ if and only if 
$\path_{T}(u_1,v_1) \cap 
 \path_{T}(u_2,v_2) = \emptyset$ 
for all $u_1, v_1 \in L_1$ and $u_2, v_2 \in L_2$, as needed.
$\QED$

%%%%%%%%%%% SHORT
%}
%%%%%%%%%%%%%%%%%% END SHORT

%%%%%%% LONG
Next we sate a useful closure property.
%%%%%%% END LONG 

%%%%%%%%%% SHORT
%Next we sate an easy and  useful closure property (proof omitted).
%%%%%%%%%% END SHORT

\begin{Lemma} \label{lem:path_edge}
If $T|L_1$ and $T|L_2$ are edge sharing and 
$L_1 \cup L_2 \subset \L(T)$, then every
edge of $e$ of $T|{L_1 \cup L_2}$ belongs to a path 
$\path_{T | L_1 \cup L_2}(u,v)$ where $u,v \in L_1$ or $u,v \in L_2$. 
\end{Lemma}

%%%%%%%%%%%%%%%%%%%%%%%% SHORT
%\ignore{
%%%%%%%%%%%%%%%%%%%%%%%%% END SHORT
\proofs
Suppose otherwise and let $e=(w,w')$ be an edge of $T | L_1 \cup L_2$ 
that does not belong to any such path. 
It follows that all the vertices of $L_1$
are on one side of that edge and all the vertices of $L_2$ on the other
side. 

Thus $\path_{T | L_1 \cup L_2}(u_1,v_1) \cap 
\path_{T | L_1 \cup L_2}(u_2,v_2) = \emptyset$ 
for all $u_1,v_1 \in L_1$ and $u_2,v_2 \in L_2$. 
This in turn implies by Lemma \ref{lem:prop_disj} that $L_1$ and $L_2$
are edge disjoint in contradiction to our assumption.
$\QED$
%%%% SHORT %%%%%%%%%%%%%%%%%%
%}
%%%%%%%%%%%%%%%%%%%%%%%%%%%%%%% END SHORT

\begin{Lemma} \label{lem:edge_close}
Suppose that $(T|L_1,T|L_2)$ are edge sharing while 
$(T|L_1,T|L_3)$ and $(T|L_2,T|L_3)$ are edge disjoint, and let 
$L = L_1 \cup L_2$. Then $(T|L,T|L_3)$ are edge disjoint.
\end{Lemma}

%%%%%%%%%%%%% SHORT 
%\ignore{
%%%%%%%%%%%%% END SHORT
\proofs
Suppose otherwise. Let $u,u' \in L$ and $v,v' \in L_3$ such that 
$\path_{T | L}(u,u') \cap \path_{T | L}(v,v') \neq \emptyset$. 
Let $e$ be an edge that
belongs to their intersection. By the previous lemma, it follows that
there exists $w,w' \in L_1$ or $w,w' \in L_2$ such that $e \in
\path_{T | L}(w,w')$. 
Now $\path_{T | L}(w,w') \cap \path_{T | L}(v,v')$ is not empty - 
in contradiction to the fact that 
$(T|L_1,T|L_3)$ and $(T|L_2,T|L_3)$ are edge disjoint. 
$\QED$
%%%%%%%%%%%%% SHORT  
%}
%%%%%%%%%%%%% END SHORT

We note that for binary trees, the notions of edge
disjointness and vertex disjointness coincide. Let $T$ be a tree and
$L_1, L_2 \subset \L(T)$. We say that $T|L_1$ and $T|L_2$ are 
{\em vertex-disjoint} if $\path_T(u_1,v_1)$ and $\path_T(u_2,v_2)$
have no vertices in common for all $u_1,v_1 \in L_1$ and $u_2,v_2 \in L_2$. 
\begin{Proposition} \label{prop:edge_vertex_disjoint}
Let $T$ be a tree and let $L_1,L_2 \subset \L(T)$. Then if $T |
L_1$ and $T | L_2$ are vertex-disjoint, they are also edge-disjoint

Let $T$ be a tree where all the internal nodes are of degree $2$ or
$3$ and let $L_1,L_2 \subset \L(T)$. Suppose furthermore that $T |
L_1$ and $T | L_2$ do not consist of a single vertex. Then if $T |
L_1$ and $T | L_2$ are edge-disjoint, they are also vertex-disjoint
\end{Proposition}

%%%%%%%%%%%%% SHORT
%\ignore{
%%%%%%%%%%%%% END SHORT
\proofs
If two paths share an edge they also share the two end points of that
edge, so the first claim follows. 

For the second claim, suppose that $T | L_1$ and $T | L_2$ are edge
disjoint but have the vertex $v$ in common. If $v$ is a leaf, then
both $T|L_1$ and $T|L_2$ share the edge adjacent to that leaf - a
contradiction. 

If $v$ is not a leaf, then there are $u_1,v_1 \in L_1$ and $u_2,v_2
\in L_2$ such that $v \in \path_T(u_1,v_1) \cap \path_T(u_2,v_2)$. But
the degree of $v$ is at most $3$, therefore the two paths 
$\path_T(u_1,v_1)$ and $\path_T(u_2,v_2)$ have non-empty edge
intersection - a contradiction. The proof follows.
$\QED$
%%%%%% SHORT %%%%%%%%%%
%}
%%%%%%%%%%%%%%%%%%%%%%% END SHORT

Edge disjoint trees naturally define a forest.
\begin{Lemma} \label{lem:forest}
Let $T$ be a binary tree.
Let $L_1,\ldots,L_{\alpha}$ be a partition of $\L(T)$ and let 
$(T_{\gamma} = T | L_{\gamma})_{\gamma = 1}^{\alpha}$ be a collection
of (pairwise) edge disjoint trees. Then the tree $T$ may be obtained
from $T_1,\ldots,T_{\alpha}$ by adding $\alpha - 1$ edges.
%(``Adding an edge'' means connecting two leaves, 
%connecting a leaf to the midpoint of an
%edge, thus splitting the edge, or connecting two edge by an edges,
%thus splitting both edges)
\end{Lemma}

%%%%%%%%%%%%% SHORT
%\ignore{
%%%%%%%%%%%%% END SHORT
\proofs
Note first that if $L = L_1 \cup L_2$ and $T|L_1,T|L_2$ are edge
disjoint then $L_1|L_2 \in \Sigma(T)$. Therefore, in this case, $T$
may be obtained from $T|L_1$ and $T|L_2$ by adding a single edge.

In the general case, define $\ell_T( T | L_{\beta}, T | L_{\gamma})$ by 
\[
\min\{\ell_T(u,v) : u \in \path_T(u',u''),\,\, v \in \path_T(v',v''),\,\, 
                  u',u'' \in L_{\beta}, \,\, v',v'' \in L_{\gamma}\}.
\]
Take $\beta \neq \gamma$ that minimize the
distance $\ell(T | L_{\beta}, T | L_{\gamma})$ among all pairs 
$(\beta,\gamma)$.

It is easy to see that $T | L_{\beta} \cup L_{\gamma}$ is edge
disjoint from $T_{\beta'}$ for $\beta' \notin \{\beta,\gamma\}$. 
The general case follows by induction.

$\QED$
%%%%%%%%%%%%% SHORT  
%}
%%%%%%%%%%%%% END SHORT

\begin{figure}[h] \begin{center} \label{fig:ex}
\input{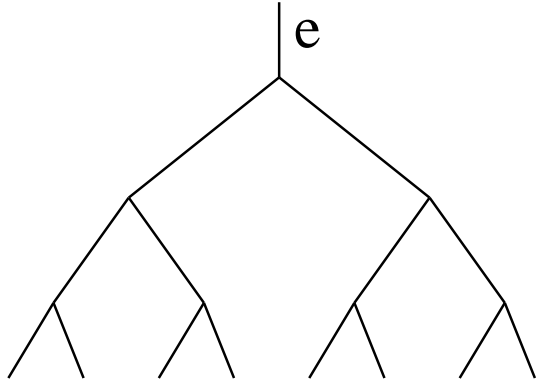}
\caption{Note that $\ell^{\ast}_T(e = (u,v),\L(T)) = 3$ while 
$\min\{\ell(u,\L(T)),\ell(v,\L(T))\} = 0$}. 
\end{center}
 \end{figure}

We say that an edge $e \in \edges(T)$ belongs to $T|L$, if there exist 
$u,v \in L$ such that $e \in \path_T(u,v)$. We say that the directed
edge $\overrightarrow{e}$ belongs to $T|L$ if the edge $e$ belongs to
$T|L$. The distance of directed edge $\overrightarrow{e} = 
\overrightarrow{(u_1,u_2)}$ to a set of
vertices $V' \subset \vert(T)$ is the minimal length $m-1$ of a {\em simple}
path $(u_1,u_2),\ldots,(u_{m-1},u_m)$ such that $u_m \in V'$. 
We
denote this distance by $\overline{\ell}_T(\overrightarrow{e},V')$.
Finally, let $\ell_T^{\ast}(e,V') = 
\max 
\left(\overline{\ell}_T(\overrightarrow{e},V'), 
          \overline{\ell}_T(\overleftarrow{e},V') \right)$
(the $\max$ is over the two orientations of $e$).
\begin{Lemma} \label{lem:iso}
Let $T$ be a binary tree and 
$L_1,\ldots,L_{\alpha}$ be a partition of $\L(T)$. 
Suppose that $(T|L_{\beta})_{\beta=1}^{\alpha}$ is a collection of
edge disjoint trees and that for all edges 
$e \in \edges(T)$ with 
$\ell^{\ast}_T(e,\L(T)) \leq r$ 
the edge $e$ belongs to one of the trees 
$T | L_{\beta}$. Then $\alpha \leq 1 +  30 \times 2^{-r} n$.
\end{Lemma}
Note that $\ell_T^{\ast}(e = (u,v),V') \geq \min\{\ell(u,V'),\ell(v,V')\}$ 
and strict inequality may hold (see Figure \ref{fig:ex}).
Thus the lemma {\em does not} follow from the fact that fractions of 
{\em vertices} at distance $r$ from the set of leaves is at most $2^{-r}$.

\proofs
Following the argument of the previous lemma, it is easy to see that
each edge we add must satisfy $\ell^{\ast}_T(e,\L(T)) > r$. 
Let $A_m$ be the set of all edges whose $\overline{\ell}$ 
distance to $\L(T)$ is at least $r$. Then $\alpha - 1 \leq |A_r|$. 
It remains to bound the size of $A_r$. 

Let $b(n,m)$ be the maximal possible size of $A_m$ among 
all binary trees on $n$ leaves. 
Note that $b(n,m) = 0$ if $n \leq 2^{m-1}$. Let $T$ be a tree on
$n$ leaves and $e$ an edge of $T$. Let $T_1$ and $T_2$ be the two
trees obtained from $T$ by removing the edge $e$. 
%(and contracting the
%degree two nodes adjacent to $e$).
Note that the number
of directed edges in $T$ of $\overline{\ell}$ distance at least $m$
from $\L(T)$ is at most five more than the sum of the number of such edges in
$T_1$ and $T_2$ (we may add at most the new edge and the four new edges
adjacent to it).

For every binary tree on $n$ leaves there exists an edge $e$ such that
removing the edge $e$ results in two tree $T_1,T_2$ such that 
$|\L(T_1)| \geq |\L(T)|/3$ and 
$|\L(T_2)| \geq |\L(T)|/3$. 
We therefore conclude that 
\[
b(n,m) \leq \max \{ b(n_1,m) + b(n_2,m) + 5 : n_1 + n_2 = n 
                    \mbox{ and } n_1 \geq n/3 \mbox{ and } n_2 \geq
                    n/3\}.
\]
It now follows by easy induction that 
\[
b(n,m) \leq \max\{30 \times 2^{-m} n - 5,0\}
\]
for all $n$. In particular,
\[ 
\alpha \leq 1 + |A_r| \leq 1 + \max\{30 \times 2^{-r} n - 5,0\} \leq 
1 + 30 \times 2^{-r} n.
\] 
as needed.
$\QED$
% For directed edges (if needed) can prove by induction - by splitting the
% tree to two.
%
%
%
 
%%%%%%%%%%%%%%%%%%%%%%%%%%%%%%%%%%%%%%%%%%%%%%%%%%%%%%%%%%%%%%%%%%%%%
%%%%%%%%%%%%%%%%%%%%%%%%%%%%%%%%%%%%%%%%%%%%%%%%%%%%%%%%%%%%%%%%%%%%%
\section{Super-trees for edge sharing trees} \label{sec:super_tree}
%%%%%%%%%%%%%%%%%%%%%%%%%%%%%%%%%%%%%%%%%%%%%%%%%%%%%%%%%%%%%%%%%%%%%
%%%%%%%%%%%%%%%%%%%%%%%%%%%%%%%%%%%%%%%%%%%%%%%%%%%%%%%%%%%%%%%%%%%%%
In this section we show how to build the super-tree of a
collection of edge-sharing trees. 

\begin{Definition} \label{def:share_collection}
Let $T$ be a binary tree and 
$L_1,L_2,\ldots,L_{\alpha} \subset \L(T)$. We say that 
$T | L_1,\ldots,T | L_{\alpha}$ are {\em edge sharing} if there is no
partition $S_1 \cup S_2$ of $\{1,\ldots,\alpha\}$ such that $T |
L_{\beta}$ and $T | L_{\gamma}$ are edge disjoint for all $\beta \in S_1$ and $\gamma \in S_2$.
\end{Definition}

From Lemma \ref{lem:edge_close} it follows that 

\begin{Proposition} \label{prop:share_collection}
Let $T$ be a binary tree and 
$L_1,L_2,\ldots,L_{\alpha} \subset \L(T)$. Then 
$T | L_1,\ldots,T | L_{\alpha}$ are edge sharing if and only if there is no
partition $S_1 \cup S_2$ of $\{1,\ldots,\alpha\}$ for which 
$T | \cup_{\beta \in S_1} L_{\beta}$ and 
$T | \cup_{\beta \in S_2} L_{\beta}$ are edge disjoint.
\end{Proposition}

In the main result of this section we prove the following
\begin{Theorem} \label{thm:super_tree}
Let $T$ be a binary tree and 
$L_1,L_2,\ldots,L_{\alpha} \subset \L(T)$ such that  
$T | L_1,\ldots,T | L_s$ are edge sharing. 
Let $L' = \cup_{\beta=1}^{\alpha} L_i$ and $T' = T | L'$.
For $1 \leq \beta \leq \alpha$, let 
\[
S(\beta) = \{\gamma : T|L_{\beta} \mbox{ and } 
            T|L_{\gamma} \mbox{ are edge sharing }\},
\]
and let $SL_{\beta} = \cup_{\gamma \in S(\beta)} L_{\gamma}$.
Then 
\begin{itemize}
\item
The tree $T'$ is determined by the trees 
$(T | SL_{\beta})_{\beta=1}^{\alpha}$.
\item
Moreover, given the trees $(T | SL_{\beta})_{\beta=1}^{\alpha}$, 
there is a polynomial time algorithm that computes $T'$.
\end{itemize}
\end{Theorem}

Theorem \ref{thm:super_tree} states that it is possible to glue together a 
collection of edge-sharing trees, given some ``local'' tree structures.

\begin{Lemma} \label{lem:eq_partition}
Assume the setting of Theorem \ref{thm:super_tree}. Let $e \in \edges(T
| SL_{\beta})$ satisfy that  
there exist $u,v \in L_{\beta}$ with 
$e \in \path_{T | SL_{\beta}}(u,v)$. Let $A_0|B_0$ be the partition defined
by $e$ on $SL_{\beta}$. Then there exists a unique partition 
$A|B \in \Sigma(T')$ such that
$A_0 = A \cap SL_{\beta}$ and $B_0 = B \cap SL_{\beta}$.

Moreover, for every edge $\tilde{e} \in \edges(T')$, there exists 
$1 \leq \beta \leq \alpha$, leaves 
$u,v \in L_{\beta}$ and $e \in \path_{T | SL_{\beta}}(u,v)$ such that the
partition of $T | SL_{\beta}$ 
defined by $e$ is given by $A \cap SL_{\beta} | B \cap SL_{\beta}$,
where $A | B$ is the partition defined by $\tilde{e}$. 
\end{Lemma}

%%%%%%%%%%%%%%%%%%%
%\ignore{
%%%%%%%%%%%%%%%%%%%
\proofs
Since $T' | SL_{\beta} = T | SL_{\beta}$, there exists a
partition $A|B$ of $L'$ corresponding to an edge of $T'$ 
which satisfies $A \cap SL_{\beta} = A_0$ and $B \cap SL_{\beta} =
B_0$. Thus, in order to prove the first claim, 
it remains to show that $A|B$ is unique. 

Write $e = (x_1,x_{k+1})$ and let $(x_1,x_2),\ldots,(x_k,x_{k+1})$ be the
path corresponding to $e$ in $T'$. 
Note that $A|B$ defines a partition satisfying 
$A_0 = A \cap SL_{\beta}$ and $B_0 = B \cap SL_{\beta}$ if and only if $A|B$
corresponds to one of the edges $(x_m,x_{m+1})$. The last claim 
follows from the fact that removing an edge of $T'$ 
that doesn't correspond to any edge in $T | SL_{\beta}$ 
induces the trivial partition on $SL_{\beta}$ and
removing an edge $\tilde{e}$ that correspond to the edge $e$ of
$T | SL_{\beta}$ corresponds to the partition defined by $e$.

Therefore, it suffices to show that $k=1$. 

Suppose that $k > 1$. Since the edge $(x_1,x_2)$ is not defined in 
$T' | SL_{\beta}$, and the collection $(T |
S_{\beta})_{\beta=1}^{\alpha}$ is edge-sharing, it follows that 
there exists $(u,v) \in  L_{\gamma}$ such that 
such that $(x_1,x_2) \in \path(u,v)$ and $\gamma \notin S(\beta)$.
But the fact that $(x_1,x_2) \in \path(u,v)$ implies that $T |
L_{\beta}$ and $T | L_{\gamma}$ are edge sharing - a contradiction.
Therefore $k=1$ and the first claim follows.

For the second claim, note that 
Lemma \ref{lem:path_edge} and Proposition
\ref{prop:share_collection} imply that for all $\tilde{e} \in \edges(T')$,
there exists $1 \leq \beta \leq \alpha$ and $u,v \in L_{\beta}$ such that 
$\tilde{e} \in \path_{T'}(u,v)$. By the previous argument, 
the edge $\tilde{e}$ corresponds to a
unique edge $e \in \edges(T | SL_{\beta})$, as needed.
$\QED$
%%%%%%%%%%%%% END OF IGNORE %%%%%%%%%%%%%%%%%%%%
%}
%%%%%%%%%%%%%%%%%%%%%%%%%%%%%%%%%%%%%

\prooft{Of Theorem \ref{thm:super_tree}}
We will show how to reconstruct $T'$ in time polynomial in $n$. Clearly, it
suffices to show how to reconstruct $\Sigma(T')$ in polynomial time. 
From Lemma \ref{lem:eq_partition} it follows that in order to find 
$\Sigma(T')$ it suffices to find for all $1 \leq \beta \leq \alpha$
and all edges $e \in \edges(T | SL_{\beta})$ which satisfy 
$e \in \path_{SL_{\beta}}(u,v)$ for $u,v \in L_{\beta}$:
\begin{itemize}
\item
The partition $A_0|B_0$ of $SL_{\beta}$ corresponding to the edge
$e$.
\item
The unique partition $A|B$ of $L'$ satisfying 
$A_0 = A \cap SL_{\beta}$ and $B_0 = B \cap SL_{\beta}$.
\end{itemize}
Given $T | SL_{\beta}$ it is trivial to find the partition $A_0|B_0$
of $SL_{\beta}$ corresponding to $e$. 
All that remains to show is how to find the partition $A|B$
corresponding to the edge $e$ in $T'$. This is the unique partition
satisfying $A_0 = A \cap L'$ and $B_0 = B \cap L'$.

For $S' \supset S(\beta)$ let 
$\tilde{L} = \cup_{\gamma \in S'} L_{\gamma}$.
Lemma \ref{lem:eq_partition} allows to identify the edge $e$ in 
$T | SL_{\beta}$ with the unique edge corresponding to $e$ in $T |
\tilde{L}$. We will use this identification below.

We now give an inductive construction of $A$ and $B$. 
Let $S_0(\beta) = S(\beta)$ and continue inductively by letting
for $d \geq 1$
\begin{itemize}
\item
$
S_d(\beta) = \{\gamma :\exists \delta \in S_{d-1}(\beta) 
               \mbox{ s.t. } (L_{\delta},L_{\gamma}) 
              \mbox{ are edge sharing }\} \setminus 
              \left( \cup_{c < d} S_{c}(\beta) \right). 
$
\item
Let $A_d = A_{d-1} \cup_{\gamma \in S_d^A(\beta)} L_{\gamma}$, where 
$\gamma \in S_d(\beta)$
belongs to $S_d^A(\beta)$ if the following condition holds:
There exists leaves
$u',v' \in A_{d-1}$ and leaves $u,v \in L_{\gamma}$ such that 
$\path(u',v') \cap \path(u,v) \neq \emptyset$.
\item
Similarly, let  
$B_d = B_{d-1} \cup_{\gamma \in S_d^B(\beta)} L_{\gamma}$, 
where $\gamma \in S_d(\beta)$
belongs to $S_d^B(\beta)$ if the following condition holds:
There exists leaves  
$u',v' \in B_{d-1}$, leaves $u,v \in L_{\gamma}$ such that 
$\path(u',v') \cap \path(u,v) \neq \emptyset$.
\end{itemize}
The above construction is repeated until $S_d(\beta) = \emptyset$.

We now prove the validity of the construction.
First, from the fact that $T|L_1,\ldots,T|L_{\alpha}$ are edge sharing, it
follows that $\{1,2,\ldots,\alpha\} = \cup_{d \geq 0} S_d(\beta)$. We
write $\bar{S}_d(\beta)$ for $\cup_{c \leq d} S_d(\beta)$. 

\begin{Claim}
For all $d \geq 0$ it holds that $S_d(\beta) = S_d^A(\beta) \cup
S_d^B(\beta)$. For all $d \geq 0$, the partition 
$A_d | B_d$ is the partition of   
the tree $T | \cup_{\gamma \in \bar{S}_d(\beta)} L_{\gamma}$
defined by the edge $e$.
\end{Claim}
\proofs
The proof is by induction on $d$. The base case $d=0$ is immediate.
For the inductive step note that under the induction hypothesis for
$d-1$, the partition $A_{d-1} | B_{d-1}$ is the partition induced by
$e$ on the tree $T | \cup_{\gamma \in \bar{S}_{d-1}(\beta)}
L_{\gamma}$. 

Let $\gamma \in S_d(\beta)$. Clearly $T | L_{\gamma}$ and 
$T | A_{d-1} \cup B_{d-1}$ share edges. On the other hand since $e$
is not an edge of $T|L_{\gamma}$ it follows that either $T | L_{\gamma}$
and $T | A_{d-1}$ share edges or $T | L_{\gamma}$ and 
$T | B_{d-1}$ share edges.  In the first case $L_{\gamma} \subset
A_d$, while in the second case $L_{\gamma} \subset B_d$. 

The claim follows.
$\QED$

The proof of the theorem follows as when the algorithm terminates, the
sets $A_d | B_d$ define the desired partition. It is clear that the
algorithm described above runs in polynomial time.
$\QED$

%%%%%%%%%%%%%%%%%%%%%%%%%%%%%%%%%%%%%%%%%%%%%%%%%%%%%%%%%%%%%%%%%%%%%
%%%%%%%%%%%%%%%%%%%%%%%%%%%%%%%%%%%%%%%%%%%%%%%%%%%%%%%%%%%%%%%%%%%%%
\section{Distorted metrics on trees} \label{sec:dist}
%%%%%%%%%%%%%%%%%%%%%%%%%%%%%%%%%%%%%%%%%%%%%%%%%%%%%%%%%%%%%%%%%%%%%
%%%%%%%%%%%%%%%%%%%%%%%%%%%%%%%%%%%%%%%%%%%%%%%%%%%%%%%%%%%%%%%%%%%%%
In this section we will prove Theorem \ref{thm:part}. We will 
assume (\ref{eq:fg}) below (i.e. $f \leq d(e) \leq g$  
for all $e \in E$). 

%%%%%%%%%%%%%%%%%%%%%%%%%%%%%%%%%
\ignore{
%%%%%%%%%%%%%%%%%%%%%%%%%%%%%%%%%
We reformulate Theorem \ref{thm:part} in terms of edge-disjoint trees.
\begin{Theorem} \label{thm:part2}
Let $T=(V,E)$ be a binary tree equipped with a metric $d$ satisfying 
(\ref{eq:fg}). Let $n = |\L(T)|$, so that $|\vert(T)| = 2 n - 2$. 
Let $\d$ be an $(\eps,M)$ distortion of $d$ and suppose
that $\eps < f/2$ and that $M > 7 \eps$. 
Then there exists a partition $\part$ of $\L(T)$ into sets 
$L_1,\ldots,L_{\alpha}$ such that 
\begin{itemize}
\item
$T|L_{\beta}$ and $T|L_{\gamma}$ 
are edge-disjoint and vertex-disjoint for $\beta \neq \gamma$.
\item
The partition size $\alpha$ is of size at most
$1 + \frac{60 n}{\sqrt{2}} \times 2^{-\frac{M-\eps}{2g}}$.
\item
The tree $T$ may be obtained from the forest $\{T|L_{\beta}\}_{\beta
  \leq \alpha}$, by adding at most $\alpha-1$ edges.
\end{itemize}

Moreover, 
\begin{itemize}
\item
the partition $(L_{\beta})_{\beta \leq \alpha}$,
\item 
the trees $(T|L_{\beta})_{\beta \leq \alpha}$ and 
\item
a function 
$\d : \cup_{\beta \leq \alpha} \edges(T|L_{\beta}) : \to \R_{+}$
satisfying $|\d(e)-d(e)| < 2 \eps$ for all $e$,
\end{itemize}
can be all found in time polynomial on $n$.
\end{Theorem}

%%%%%%%%%%%%%%%%%%%%% END IGNORE %%%%%%%%%%
}
%%%%%%%%%%%%%%%%%%%%%%%%%%%%%%%%%%%%%%%

It is helpful to define ``balls'' with respect to $d$ and $\d$ as follows.
\begin{equation} \label{eq:ball1}
B_L(v,r) = \{w \in \L(T) : d(v,w) \leq r\},  B_V(v,r) = \{w \in \vert(T) : d(v,w) \leq r\}.
\end{equation} 
We similarly define $\B_L(v,r)$ and $\B_V(v,r)$ with $\d$ instead of
$d$. 

We omit the proof of the following easy Lemma. 

\begin{Lemma} \label{lem:Cprop}
Let $\d$ be an $(\eps,M)$ distortion of $d$ and 
let $r < M$. Consider the sets 
$L_{\beta} = \B_L(v_{\beta},r)$ for $\beta=0,\ldots,\alpha$.
Suppose that $T | L_0$ and $T | L_{\beta}$ are edge sharing for 
$\beta=1,\ldots,\alpha$. Then for all 
$u,u' \in \vert(T | \cup_{\beta = 0}^{\alpha} L_{\alpha})$ it holds that 
\begin{equation} \label{eq:d3}
d(u,v) < 6 r + 6 \eps. 
\end{equation}
Similarly, for all 
$u,u' \in \cup_{\beta = 0}^{\alpha} L_{\alpha}$ it holds that 
\begin{equation} \label{eq:d4}
\d(u,v) < 6 r + 7 \eps,
\end{equation}
if $\d(u,v) < \infty$. 
\end{Lemma}
%%%%%%%%%%%%%%%%%%%%%%%%%%%%%%%%%%% SHORT 
%\ignore{
%%%%%%%%%%%%%%%%%%%%%%%%%%%%%%%%%%% END SHORT
\proofs
Equation (\ref{eq:d4}) follows immediately from (\ref{eq:d3}).
Let $T' = T | \cup_{\beta = 0}^{\alpha} L_{\beta}$. 
Note that if $\beta \leq \alpha$, the leaves $u_1,u_2$ belong to
$L_{\beta}$ and $u$ is a vertex on the path connecting $u_1$ and
$u_2$, then 
\begin{equation} \label{eq:stupid_radius}
d(u,v_{\beta}) \leq \sup_{w \in L_{\beta}} d(w,v_{\beta}) < 
\sup_{w \in L_{\beta}} \d(w,v_{\beta}) + \eps \leq r + \eps. 
\end{equation}
(the first inequality follows from the fact that if $v$ is a vertex in a
tree with a path metric, then the vertex furthest away from $v$ is a leaf).

Now let $u,u' \in \vert(T')$. Since the trees $(T | L_{\beta})_{\beta
  \leq \alpha}$ are edge
sharing, it follows $u$ belongs to a path connecting two points in
  $L_{\gamma}$ and $u'$ belongs to a path connecting two point in 
$L_{\gamma'}$, where $0 \leq \gamma, \gamma' \leq \beta$.

Let $w$ be a vertex that belongs to an
  edge $e$ such that 
$e \in \path(w_{\gamma,1} , w_{\gamma,2}) \cap \path(w_{0,1} , {w}_{0,2})$
where $w_{\gamma,1},w_{\gamma,2} \in L_{\gamma}$ and  
$w_{0,1},w_{0,2} \in L_{0}$. Define $w'$ similarly.
Then by (\ref{eq:stupid_radius})
\begin{eqnarray*}
d(u,u') &\leq& d(u,v_{\gamma}) + d(v_{\gamma},w) + 
             d(w,v_{0}) + d(v_{0},w') + 
             d(w',v_{\gamma'}) + d(v_{\gamma'},u') \\
        &\leq& 6 (r + \eps),
\end{eqnarray*}
as needed.   
$\QED$ 
%%%%%%%%%%%%%%%%%%%%%%%%%%%%%%%%%%%%%% SHORT 
%}
%%%%%%%%%%%%%%%%%%%%%%%%%%%%%%%%%%%%%  END SHORT

\prooft{of Theorem \ref{thm:part}}
Let 
\[
\begin{array}{lll}
r = \frac{M - 7 \eps}{6},& L = (v_{\beta})_{\beta=1}^n,&
L_{\beta} = \B(v_{\beta},r).
\end{array}
\]
Define a graph $G$ on the set of vertices $\{1,\ldots,n\}$ where the edge 
$(\beta,\gamma)$ is present if and only if $L_{\beta}$ and
$L_{\gamma}$ are edge-sharing. Let $C_1,\ldots,C_{\alpha}$ be the partition
of $\{1,\ldots,n\}$ to connected components in $G$. 
It is easy to see that the graph $G$ can be computed in polynomial
time. Indeed by Lemma \ref{lem:Cprop} it follows that if $L_{\beta}$
and $L_{\gamma}$ share edges then for all $u \in L_{\beta}$ and $v \in
L_{\gamma}$ it holds that $\d(u,v) \leq M$. For sets
$L_{\beta},L_{\gamma}$ for which this condition holds, we may easily
reconstruct the tree $T | L_{\beta} \cup L_{\gamma}$ using the
$4$-point method. We can then check if $T | L_{\beta}$ and $T |
L_{\gamma}$ are edge sharing in $T | L_{\beta} \cup L_{\gamma}$.

By proposition (\ref{prop:share_collection}) it follows that if 
$\sigma \neq \tau$ then the trees 
$T | \cup_{\beta \in C_{\sigma}} L_{\beta}$ and
$T | \cup_{\beta \in C_{\tau}} L_{\beta}$ are edge disjoint. 
Moreover, if $1 \leq \eta \leq \alpha$, then the collection of trees 
$(T | L_{\beta})_{\beta \in C_{\eta}}$ are edge sharing. 

Note that in the notation of Theorem \ref{thm:super_tree} for all
$\beta$ and all $u,v \in SL_{\beta}$ it holds that 
$\d(u,v) \leq M$ - by Lemma \ref{lem:Cprop}. 

It follows that the trees $T | SL_{\beta}$ may be easily recovered by
the $4$-point method. Moreover, for every edge $e \in \edges(T |
SL_{\beta})$ we may recover $\d(e)$ satisfying $|d(e) - \d(e)| < 2
\eps$. 

We now use Theorem \ref{thm:super_tree} in order to recover
the trees $T | \cup_{\beta \in C_{\sigma}} L_{\beta}$ for all
$\sigma$. Moreover, we may recover $\d(e)$ satisfying $|d(e) - \d(e)| < 2
\eps$. 

It remains to bound the number of trees $\alpha$ using Lemma
\ref{lem:iso}. Note that if $l^{\ast}_T((u,v),\L(T)) \leq p$, then
there is a path of $p-1$ edges starting at $u$, avoiding $v$ and
ending at $\L(T)$ at a node denoted $u'$. 
Similarly, there is a path of $p-1$ edges starting at $v$ avoiding $u$
and ending at $\L(T)$ at a node denoted $v'$. Note that 
$d(u',v') \leq (2p-1)g$ and therefore $\d(u',v') \leq (2p-1)g +
\eps$. Thus if $p \leq \frac{1}{2} + \frac{M - \eps}{2g}$, then 
$\d(u',v') \leq M$ and therefore the edge $e$ belongs to one of the
trees $T | \cup_{\beta \in C_{\sigma}} L_{\beta}$.  
It follows from Lemma \ref{lem:iso} that
\[
\alpha \leq 
1 + 30 n \times 2^{- \lfloor \frac{1}{2} + \frac{M - \eps}{2g} \rfloor} \leq 
1 + \frac{60 n}{\sqrt{2}} \times 2^{-\frac{M-\eps}{2g}}.
\]  
The theorem follows. 
$\QED$

%%%%%%%%%%%%%%%%%%%%%%%%%%%%%%%%%%%%%%%%%%%%%%%%%%%%%%%%%%%%%%%%%%%%%%%%%%%%%%%%%%%%%
%\bibliographystyle{abbrv}
%\begin{thebibliography}{99}
%%%%%%%%%%%%%%%%%%%%%%%%%%%%%%%%%%%%%%%%%%%%%%%%%%%%%%%%%%%%%%%%%%%%%%%%%%%%%%%%%%%%

\bibliographystyle{plain}
\bibliography{all}
 
%\end{document}

\appendix
\section{Large deviations for the log-det distance}
%Note that if $\det(M(e)) \in \pm 1$, then $M(e)$ is a permutation
%matrix and therefore there is no mutation along $e$. If $\det(M(e)) =
%0$, then $M$ ``losses information'' - a case which we also like to
%avoid. Finally note by relabelling some of the states if needed, we
%may assume that $\det(M(e)) > 0$ for all $e$. In fact we will assume
%more:
%\begin{itemize}
%\item 
%There exists $g > f > 0$ such that $e^{-2g} < \det(M(e)) < e^{-2f}$
%for all $e \in \edges(T)$.
%\item
%There exists a positive number $\pi_{\min}$ such that for all nodes
%$v \in \vert(T)$ and all letters $a \in \A$, the probability that the
%letter at $v$ is $a$ is at least $\pi_{\min}$.
%\end{itemize}

In this section, we prove Proposition \ref{prop:large_deviations}.
We will use the following large deviation result from \cite{ErStSzWa:99b}. 
\begin{Lemma}[\cite{ErStSzWa:99b}] \label{lem:essw2}
For every two vertices $u,v \in \L(T)$ and $a,b \in \A$, 
let $F_{a,b}$ be the probability that node $u$ has letter $a$ and node
$w$ has letter $b$. Let $-d(u,v) = \log \det (F_{i,j})$. Then $d(u,v)$
is a path metric on the tree satisfying 
\begin{itemize} 
\item 
$d(e) \geq f$ for all edges $e$ of the tree.
\item
For all $u,v \in \L(T)$:
\begin{eqnarray*}
d(u,v) &\leq& -|\A| \log \pi_{\min} - 
             \sum_{e \in \path(u,v)} \log \det(M(e)) \\ &\leq&  
- \|A| \log \pi_{\min} + 2g |\path(u,v)|.
\end{eqnarray*}
\end{itemize}
Moreover for $u,v \in \L(T)$ let $\F_{a,b}$ be the empirical
distribution of having $a$ at $u$ and $v$ in $b$ in a collection of
$k$ samples. Let $\d(u,v) = -\log \det(\F_{a,b})$ if 
$\det(\F_{a,b}) >  0$ and $\d(u,v) = \infty$ otherwise. Then there
exists positive constants $c_1$ and $c_2$ such that 
\begin{equation} \label{eq:large_deviations}
\P[|e^{-d(u,v)} - e^{-\d(u,v)}| \geq t] \leq 2 \exp 
\left( -c_1 k \left(t - \frac{c_2}{k} \right)_+^2 \right),  
\end{equation}
where $(a)_+ = \max\{0,a\}$. 
\end{Lemma}
For the proof see \cite[Section 7]{ErStSzWa:99b}. Equation
(\ref{eq:large_deviations}) here is equation (49) in \cite{ErStSzWa:99b} up
to change of notation. Taking $t = e^{-M-\eps}(1 - e^{-2\eps})$ in 
(\ref{eq:large_deviations}) we see that if 
$|e^{-d(u,v)}-e^{-\d(u,v)}| \leq t$
and either $d(u,v) \leq M+\eps$ or $\d(u,v) \leq M+\eps$, then 
$|d(u,v)-\d(u,v)| < \eps$. So if 
$|e^{-d(u,v)}-e^{-\d(u,v)}| \leq t$ for all $u$ and $v$, then $\d$ is
an $(\eps,M)$ distortion of $d$. 

Taking $k$ that satisfies (\ref{eq:2nd_dev}) we obtain that the error
is at most 
\[
n^2 \exp \left( -c_1 k \left(t - \frac{c_2}{k} \right)_+^2 \right) \leq 
n^{2-r},
\]
if $c$ (in (\ref{eq:2nd_dev})) is sufficiently large. 
The proof of Proposition \ref{prop:large_deviations} follows.

\section{Lower bounds}
In this section we prove tightness of both the distorted metric result 
and the phylogenetic reconstruction result. 

The tightness of the metric result follow easily by considering the 
$r$-level $3$-regular tree with the metric $d$ that assigns length $d$ to all 
edges of the tree. We let $\d(u,v) = d(u,v)$ if $d(u,v) \leq M$ and 
$\d(u,v) = \infty$ otherwise. Then $\d$ is a $(0,\infty)$ distortion of $d$.

Define the relation $u \sim v$ if $d(u,v) \leq M$. It is easy to see 
that $\sim$ is an equivalence relation. There are
$n 2^{-\lfloor M/2g \rfloor}$ equivalence classes for this relation. 
It is easy to reconstruct the tree on each class, but since for $u,v$ 
which belong to different classes, $\d(u,v) = \infty$, it is impossible 
to reconstruct any more. This prove the tightness of the number of trees in 
Theorem \ref{thm:part} up to a multiplicative constant.

A similar construction yields the analogous sampling 
complexity lower bound for phylogenetic trees. We fix the model to be 
the CFN model where the length of each edge is $g$. Thus the mutation matrices 
are given by 
$M(e) = \left( \begin{matrix} e^{-g} + (1 - e^{-g})/2 & (1-e^{-g})/2 \\ 
                                (1-e^{-g})/2 & e^{-g} + (1-e^{-g})/2 
\end{matrix} \right)$. 

Thus for each edge $(u,v)$, the state of $u$ is copied to $v$ 
with probability $e^{-g}$. Otherwise, an independent uniform state is chosen. 

Following the arguments of \cite{Mossel:03b} 
implies that if $v$ is a vertex at 
$\ell$-distance $s$ from the set of leaves, then the character value at the 
leaves below $v$ is independent of the character at $v$ with probability at 
least $1 - 2^s e^{-g s}$. Thus the character at the leaves 
is independent from all nodes at level $s$ 
with probability at least $1 - 3*2^{r-s-1} 2^2 e^{-g s} = 
1 - n e^{-g s}$. The probability that the former event will occur for $k$ characters is at least $p_k = 1 - k n e^{-g s}$. 

Let assume further that the phylogenetic tree on each of the equivalence 
classes of the relation $\sim$ defined by $u \sim v$ if $d(u,v) \leq 2 g s$ is given.
Then with probability $p_k$, there is no non-trivial information about the 
ancestral relationship except that given by the given $n*2^{-s}$ trees. 

Note furthermore that $p_k \geq 1 - \delta$ if $k \leq \delta e^{g s}/n = 
\delta e^M/n$. This proves the tightness of condition (\ref{eq:2nd_dev}) in 
Theorem \ref{thm:phyl} up to a factor $2$ in the exponent and a multiplicative 
$O(n)$ factor.

\section{Variants of the method}
We briefly sketch a few variants of the method which may be practical 
advantages over the method analyzed here.

\subsection{Checking if two balls define tree disjoint trees}
The first stage of the algorithm consists of checking if two balls $\B_L(v,r)$ and $\B_L(u,r)$ define two 
edge-sharing trees or two edge-disjoint trees. Most of the work at this stage is devoted to couples 
of trees that are edge-disjoint. In fact for most such pairs it would hold that $\d(u,v) > M$ which implies 
automatically edge-disjointness without additional computation. Thus the efficient way of computing the graph $G$ is by first checking for each $u$ and $v$ if $\d(u,v) > M$. Otherwise, we perform the test described in the proof of Theorem \ref{thm:part}. 

\subsection{Building supertrees for edge disjoint trees}
A lot of computational effort is devoted to building super-trees from collection of edge-sharing trees. 
There are many variants that work here. Instead of the method described in the paper, we can use 
quartets method as in \cite{ErStSzWa:99a}. 
Similarly to \cite{ErStSzWa:99a} one can prove that given a collection of edge sharing trees $T_1,\ldots,T_{\alpha}$, 
their super-tree is in fact defined by all quartets belonging to the 
trees $T | SL_{\beta}$ for $1 \leq \beta \leq \alpha$ via the dyadic closure operator. This may lead to a 
computationally more efficient algorithm than ours. 

%\subsection{What to do with short distances?}
%In order to take care of short distances - we will define an
%equivalece relation between points, where two points are in the same
%equivalence claess if all edges connecting them are shorter than some
%$\eps$. We will be interested in finding a tree on these collapsed
%clusters (in practice the egdes my be between $\eps$ and $2 \eps$). 
%This can be done if for every two cluster of vertices we define the
%distance between them as the minimal distance between any two points
%in a cluster. Does it work out?

%\subsection{Forests of forests}
%Another important extenstion is given when we are given a forest to
%start with. We are allowed to compute all distance between any pair of
%internal vertices in the forest. Again we shuold be able to get a
%refinement of the original forest.

%%%%%%%%%%%%%%%%%%%%%%%%%%%%%%%%%%%%%%%%%%%%%%%%%%%%%%%%%%%%%%%%%%%%%
%%%%%%%%%%%%%%%%%%%%%%%%%%%%%%%%%%%%%%%%%%%%%%%%%%%%%%%%%%%%%%%%%%%%%
\end{document}